\documentclass[11pt,a4paper,twoside]{amsart}

\usepackage{amsfonts, indentfirst}
\usepackage{ifthen,amssymb}
\usepackage[all]{xypic}
\usepackage{young}
\usepackage{graphicx}
\usepackage{hyperref}

\usepackage{tikz}
\usetikzlibrary{automata}

\newtheorem{prop}{Proposition}[section]
\newtheorem{theorem}[prop]{Theorem}
\newtheorem{lemma}[prop]{Lemma}
\newtheorem{example}{Example}
\newtheorem{defi}[prop]{Definition}
\newtheorem{coro}[prop]{Corollary}
\newtheorem{remark}[prop]{Remark}

\newcommand{\cqd}{\hfill$\Box$}

\renewcommand{\dim}[0]{\operatorname{dim}}

\newcommand{\PG}[0]{\operatorname{PG}}
\newcommand{\PGL}[0]{\operatorname{PGL}}
\newcommand{\GL}[0]{\operatorname{GL}}

\title[Hypergraph encoding set systems and their linear representations]{Hypergraph encoding set systems and their linear representations}
 \author[Cristina Mart{\'\i}nez]{Cristina Mart{\'\i}nez}
 \author[Alberto Besana]{Alberto Besana}
 

 \email{cristina.martinez@inv.uam.es}
\email{alberto.besana@gmail.com}
\begin{document}
\maketitle
\begin{abstract}
We study $t$-designs of parameters $(n,k,\lambda)$ over finite fields as group divisible designs and set systems admitting a transitive action of a linear group encoded in an hypergraph $G$ whose vertex set of size $n$ is partitioned into sets of size $k$ in such a way that every $t$-subset is contained in at least $\lambda$ subsets of $G$.%relating the problem to the representation theory of the general linear group $\GL(n,\mathbb{F}_{q})$ and the constructions of AG codes over finite fields. 
%We outline a method to construct large sets of $t-$designs of arbitrary size. 
We relate the problem to the representation theory of the general linear group $\GL(n,\mathbb{F}_{q})$ and the constructions of AG codes over finite fields. 
As a byproduct we construct a RS code based encryption scheme.
%
%codes and $t-$designs with an action of an abelian $p-$group.
%construct families of cyclic codes by considering sets of roots of the polynomial $x^{n}-1$ over its splitting field.

\end{abstract}

\section{Introduction}
Let $q$ be a power of a prime number and $\mathbb{F}_{q}$ the finite field with $q$ elements.  Finite fields have a remarkable property that finite dimensional vector spaces over them are naturally endowed with a canonical and compatible field structure. %In particular we have $\mathbb{F}_{q^{n}}=\mathbb{F}_{q}(\alpha)$, for some element $\alpha \in \mathbb{F}_{q}$.% and the set $\{1,\alpha, \alpha^{2},\ldots, \alpha^{m-1}\}$ forms a basis of $\mathbb{F}_{q^{m}}$ over the prime field $\mathbb{F}_{p}$.

%The encoding of an information word into a $k-$dimensional subspace is usually known as coding for errors and erasures in random network coding, \cite{KK}. %More precisely,
 %Namely, 
 We refer to network coding as the best way to disseminate information over  a network.
 %Let $V$ be an $N-$dimensional vector space over $\mathbb{F}_{q}$, a code for an operator channel with ambient space $V$ is simply a nonempty collection of subspaces of $V$. The collection of subspaces is a code for error correcting errors that happen to send data through an operator channel. %or by $\PG^{r}(n,q)$. 
A network is usually represented by a directed multigraph with error free unit capacity edges. There are several source nodes and several destination nodes and data is transferred over a network using packets, where a packet is just an $m$-length vector over a finite field $F_{q}$. %Given integers $n,k,$ and $d$, an $[n,k,d]_{q}-$code over  $\mathbb{F}_{q}$ is a subspace $C$ of $(\mathbb{F}_{q})^{n}$ of dimension $k$, such that every nonzero codeword $0\neq \alpha \in C$ satisfies $wt\,(\alpha)\geq d $, where the weight of a vector is the number of non-zero components. The minimum distance is equal to the smallest of the weights of the nonzero codewords. %$\mathbb{F}_{q}$. 
The network nodes exchange messages being represented as a matrix. % and coset coding allows to change an incoming message so that an outgoing message looks very differently.
It is convenient to describe a coding process in terms of operations in the extended field $\mathbb{F}_{q}^{m}$.
More generally, if $A$ is a commutative ring with identity, a linear code $C$ of length $n$ over $A$ is an $A$-submodule of $A^{n}$.  

Let $V$ be an $n+1$ dimensional vector space over the field $\mathbb{F}_{q}$, we denote by %${\rm{PG}}(n,q)$ or 
${\rm{\mathbb{P}}}(V)$ the $n$-dimensional projective space over it. The set of all subspaces of dimension $r$ is the Grassmannian $\mathcal{G}_{r,n}({\mathbb{F}_{q}})$ of $r$-dimensional subspaces in $(\mathbb{F}_{q})^{n}$. A subspace code is a constant dimension code (CDC), that is a subset of the Grassmannian.%Any subset of the Grassmannian is a constant dimens

In general, for any integers $n, r$ with $n\geq r \geq 0$,  we call $\phi(r;n,q):=|\PG^{r}(n,q)|$, the number of $r$ dimensional subspaces of an $n$ dimensional subspace over $\mathbb{F}_{q}$. It is  the number of ways of choosing $r+1$ linearly independent points in $\PG(n,q)$ divided by the number of ways of choosing such a set of points in a particular $r$-space. It is given by the $q$-ary binomial coefficient,

 $$\phi(q;n,r)=\left[
\begin{matrix}
 n \\
 r\\
\end{matrix}
\right]_{q}=\frac{(q^{n+1}-1)(q^{n+1}-q)\ldots (q^{n+1}-q^{r})}{(q^{r+1}-1)(q^{r+1}-q)\ldots (q^{r+1}-q^{r})}.$$

%There is a right action of 
The general linear group $\GL(n,\mathbb{F}_{q})$ acts transitively on $\mathcal{G}_{k,n}(\mathbb{F}_{q})$:

\begin{eqnarray}\label{action}
 \mathcal{G}_{k,n}(\mathbb{F}_{q}) \times \GL(n,\mathbb{F}_{q}) & \rightarrow & \mathcal{G}_{k,n}(\mathbb{F}_{q}) \\
 (\mathcal{U},A) & \rightarrow & \mathcal{U}A.
\end{eqnarray}
Observe that the action is defined independent of the choice of the representation matrix $\mathcal{U}\in \mathbb{F}_{q}^{k\times n}$. %This is known as coset coding which allows to change an incoming message so that an outgoing message looks quite differently. In particular, incoming and outgoing messages are statistically independent.

In order to classify the orbits of $\mathcal{G}_{k,n}(\mathbb{F}_{q})$ by the action of the general linear group $\GL(n,\mathbb{F}_{q})$ we need to classify all the conjugacy classes of subgroups in $\GL(n,\mathbb{F}_{q})$. %In \cite{BM1} we studied cyclic coverings of the projective line that correspond to orbits defined by a cyclic subgroup, that is a subgroup in $\GL(n,\mathbb{F}_{q})$ containing a cyclic subgroup $\mathbb{Z}_{p}$ for some prime number.

%An incidence structure is an ordered triple $\mathcal{D}=(\mathcal{P},\mathcal{B}, \mathcal{I})$ where $\mathcal{P}$ and $\mathcal{B}$ are non-empty disjoint sets and $\mathcal{I} \subseteq \mathcal{P}\times \mathcal{B}$. The elements of the set $\mathcal{P}$ are called points, the elements ofthe set $\mathcal{B}$ are called blocks and $\mathcal{I}$ is called an incidence structure.
A group divisible design (GDD) is an incidence structure $(X,\mathcal{G},\mathcal{B})$ where $X$ is a set of points, $\mathcal{G}$ is a partition of $X$ into groups, and $\mathcal{B}$ is a collection of subsets of $X$ called blocks such that any pair of distinct points from $X$ occurs either in some group or in exactly one block, but not both, see \cite{GW}. %Whenever the size of the groups is exactly 1 or each $t-$subspace of $X$ is contained in exactly $\lambda$ blocks of $\mathcal{B}$ the GDD is a block $t-$design.
Recently, $q$-ary designs (designs over finite fields) gained a lot of attention because of its applications for error-correcting in networks, and secret sharing scheme, a way for sharing a secret data among a group of participants so that only specific subsets (which are called qualified subsets) are able to recover the secret by combining their shares.

In this paper, we give an explicit method to construct large sets of $t$-designs over finite fields of given parameters $n,q, s$, by deriving ordered basis of $(\mathbb{F}_{q})^{n}$.
%We pretend to study network coding relating the problem to the representation theory of the general linear group $\GL(n,\mathbb{F}_{q})$, where $q$ is a power of a prime number $p$. 

\subsection*{Notation}
For $d$ a positive integer,  $\alpha=(\alpha_{1},\ldots, \alpha_{m})$ is a partition of $d$ into $m$ parts if the $\alpha_{i}$  are positive and decreasing. We will denote as $\mathcal{P}(d)$, the set of all partitions of $d$. We set $l(\alpha)=m$ for the length of $\alpha$, that is the number of cycles in $\alpha$, and $l_{i}$ for the length of $\alpha_{i}$. The notation $(a_{1},\ldots, a_{k})$ stands for a permutation in $S_{d}$ that sends $a_{i}$ to $a_{i+1}$. 
We write $\PGL(2,k)=\GL(2,k)/k^{*}$, and elements of $PGL(2,k)$ will be represented by equivalence classes of matrices 
$\left(\begin{array}{ll} a &  b \\
c & d \end{array}\right)$, with $ad-bc\neq 0$. 
In the sequel, %a $q-$ary constant weight code of length $n$, distance $d$ and weight $w$ will be denoted as an $[n,d,w]_{q}$ code. 
an $[n,k]_{q}$-code $C$ is a $k-$dimensional subspace of $(\mathbb{F}_{q})^{n}.$

\section{$t-(n,k,\lambda;q)$-designs}
An incidence structure with $v$ points, $b$ blocks and constant block size $k$ in which every point appears in exactly $r$ blocks is a GDD with parameters $(v,b,r,k,\lambda_{1},\lambda_{2},m,n)$ whenever the point set can be partitioned into $m$ classes of size $n$, such that two points from the same class appear together in exactly $\lambda_{1}$ blocks, and two points from different classes appear together in exactly $\lambda_{2}$ blocks, \cite{CMR}.% Then the parameters satisfy the following relations:
%$$v=mn\ \ bk=vr,\ \ (n-1)\,\lambda_{1}+n\,(m-1)\lambda_{2}=r,(k-1), \ \ rk\geq v\lambda_{2}.$$
\begin{example}\label{3-design} Consider the GDD given by the incidence structure $(\mathcal{P},\mathcal{B},I)$ on a 3-dimensional vector space $V$ over the finite field $\mathbb{F}_{p}$, where $\mathcal{P}$ is a set of $v$ smooth, reduced points in $V$, $\mathcal{B}$ is a set whose elements are triples of points $(p,q,r)\in \mathcal{P}\times \mathcal{P}\times \mathcal{P}$ defined by the condition $(p,q,r)\in \mathcal{B}$ if either $p+q+r$ is the full intersection cycle of the projective line with a $\mathbb{F}_{p}$-line $l\subset \mathbb{P}(V)(\mathbb{F}_{p})$ with the right multiplicities, or else if there exists a $\mathbb{F}_{p}-$line $l\subset V,$ such that $p,q,r \in l$, then the triple is called a plane section. 
%Two points $x=(x_{0},x_{1},x_{2})$ and $y=(y_{0},y_{1},y_{2}) \in (\mathbb{F}_{p})^{3}$ are equivalent if and only if $y_{i}=t\,x_{i}$ for $i\in \{0,1,2\}$,  for some $t\in \mathbb{F}^{*}_{p}$.
%Therefore 
The number of points in the projective plane $PG(2,p)$ is $\frac{p^{3}-1}{p-1}=p^{2}+p+1$ and dually there are $p^{2}+p+1$ lines in $PG(2,p)$. %The dual $x^{*}$ of a point in $PG(2,p)$ is a line.

%The underline incidence relation is expressed in terms of collinearity of triples of points. Namely, 
There are two types of GD designs on $V$:
\begin{enumerate}
\item For any $(p,q)\in \mathcal{P}^{2}(V^{*}),$ there exists an $r\in \mathcal{P}(S^{d}V^{*})$ such that $(p,q,r)\in l$, where $S^{d}V$ is the $d^{th}$ symmetric power of the dual of $V$. The triple $(p,q,r)$ is strictly collinear if $r$ is unique with this property, and $p,q, r$ are pairwise distinct. The subset of strictly collinear triples is a symmetric ternary relation.
\item  Assume that $p\neq q$ and there are two distinct points $r_{1}, r_{2} \in \mathcal{P}$ with $(p,q,r_{1})\in \mathcal{B}$ and $(p,q,r_{2})\in \mathcal{B}$. Denote by $l=l(p,q)$ the set of all such points, then $l^{3}\in \mathcal{B}$, that is, any triple $(r_{1},r_{2},r_{3})$ of points in $l$ is collinear. Such sets are called lines in $\mathcal{B}$.
\end{enumerate}
\end{example}

If $X$ is the finite field $\mathbb{F}_{q}^{n}$, we consider the GDD where the set points correspond to vectors of $(\mathbb{F}_{q})^{n}$ as $n$-dimensional vector space over $\mathbb{F}_{q}$ and the block set $\mathcal{B}$ is a collection of $i-$subspaces $K\subseteq \mathbb{F}_{q}^{n}$ which geometrically correspond to points in the Grassmannian  $\mathcal{G}_{n,i}(\mathbb{F}_{q})$. They live in a natural way as subspaces of the vector space $(\mathbb{F}_{q})^{n}$. More generally, assuming that $\dim K_{j}=j$ for $j=1,\ldots, n$, the sequence of nested subspaces
$$\{0\}\subset K_{1}\subset K_{2}\subset \ldots \subset K_{n}=\mathbb{F}_{q}^{n},$$ live in the whole lattice of subspaces of the vector space $(\mathbb{F}_{q})^{n}$.

If each $t$-subspace of $X$ is contained in exactly $\lambda$ blocks of $\mathcal{B}$ then it is called  a $t$-$(n,k,\lambda;q)$ design. %if lambda=1 it is called a steiner triple
A permutation matrix $\sigma \in \GL(n,q)$ acts on the Grassmannian by multiplication on the right of the corresponding representation matrix. In particular $\sigma$ is an automorphism of the design $\mathcal{D}=(X,\mathcal{B})$ if and only if $\sigma$ leaves the Grassmannian invariant, that is $\mathcal{B}^{\sigma}=\mathcal{B}$. In particular, we are interested in understanding the orbits by the action of any permutation matrix of $\GL(n,q)$ and moreover of any subgroup $G$ contained in $\GL(n,q)$. Further, it is possible to count the orbits of the action in several cases and these correspond to blocks of the design satisfying certain geometrical properties. 

 % add reference

\begin{defi} 
%\begin{enumerate}
%\item 
Let $\alpha\in \mathbb{F}_{q^{n}}$ be a generator of the underlying vector space over $\mathbb{F}_{q}$. Then an $r$-dimensional $W$  subspace is $\alpha$-splitting if $\alpha^{i}W=W$ is invariant under the action of any element $\alpha^{i}$ in the Galois group of the extension $\mathbb{F}_{q}\hookrightarrow \mathbb{F}_{q}(\alpha)$. 
%\item 
More precisely, given any $\mathbb{F}_{q}$-linear endomorphism $T:\,\mathbb{F}_{q^{n}}\rightarrow \mathbb{F}_{q^{n}}$, an $r$-dimensional subspace $W$ is $T$-splitting if $\mathbb{F}_{q^{n}}=W\oplus T(W)\oplus \cdots \oplus T^{n-1}(W),$ where $T^{j}$ denotes the $j$-fold composite of $T$ with itself.
\end{defi}
%\item 
\begin{defi}
Let $T$ be the standard shift operator on $\mathbb{F}_{q}^{n}$, a linear code $C$ is said to be quasi-cyclic of index $l$ or $l-$quasi-cyclic if and only if is invariant under $T^{l}$. If $l=1$, it is just a cyclic code. The quantity $m:=n/l$ is called the co-index of $C$. Namely, if we view a codeword $(c_{0},c_{1},\ldots, c_{n-1})$ of $C$ as a polynomial $c_{0}+c_{1}x+\ldots+c_{n-1}x^{n-1}\in \mathbb{F}_{q}[x]$, then $T(c(x))=x\cdot c(x) \ {\rm{mod}}\, (x^{n}-1)$. In particular a cyclic code $C\subset (\mathbb{F}_{q})^{n}$ is identified with an ideal in the ring $\mathbb{F}_{q}[x]/(x^{n}-1)$ generated by a polynomial $g(x)$ which divides $(x^{n}-1)$.
%\end{enumerate}
\end{defi}

\begin{example}
For $n=6$, we consider the factorization into prime factors $r=2, s=3$,  and let $\alpha$ be a primitive element of $\mathbb{F}_{q^{6}}$, then there are as many 2-subspaces of $(\mathbb{F}_{q})^{6}$ as $|\GL(2,\mathbb{F}_{q})|=(q^{3}-1)(q^{3}-q)(q^{3}-q^{2})$. There are as many 3-subspaces of $(\mathbb{F}_{q})^{6}$  as $|\GL(2,\mathbb{F}_{q})|=(q^{2}-1)(q^{2}-q)$.\,%and $\lambda=$
For $q=2$, we get a 2-(6,3,3) design. %For $q=4$, we get a 2-(24,3,3) design, a 3-(24,3,1) design, a 3-(24,4,3) design, a 3-(24,8,21) design, a 3-(24,6,10) and a 3-(24,12, 60) design.
\end{example}

One can study the orbits of $\mathcal{G}_{k,n}(\mathbb{F}_{q})$ by the action of any subgroup in the general linear group $\GL(n,\mathbb{F}_{q})$. For example we can study the orbit of any triangle group: the Klein group $\mathbb{Z}_{2}\times \mathbb{Z}_{2}$, the dihedral group, the alternated groups $A_{4}$ and $A_{5}$ or the symmetric group $S_{n}$. We study the special case of the Grassmannian $\mathcal{G}_{2,4}(\mathbb{F}_{q})$ of lines in a 3-dimensional projective space.

%\begin{prop}\label{PropSteiner} %This proposition relates the representation theory of the general linear group to algebraic codes (over finite fields)
%Invariant subgrassmannians in $\mathcal{G}_{k,n}(\mathbb{F}_{2}^{m})$ by the action of any triangle group are $t$-designs.
%\end{prop}
%{\it Proof.} Triangle groups are reflection groups which admit a presentation
%$${\rm{T}}(r,s,t)=\langle x, y, x:\, x^{r}=y^{m}=z^{w}=xyz=1\rangle,$$ with $r, m, w$ integer numbers such that $\frac{1}{r}+\frac{1}{m}+\frac{1}{w}<1$.
%They are finite subgroups of ${\rm{PGL}}(2,q)$ and they are known to be isomorphic to either the Klein group $\mathbb{Z}_{2}\times \mathbb{Z}_{2}$, the dihedral group $D_{m}$ of order $2m, m\geq 2$, the alternating group $A_{4}$, the alternating group $A_{5}$, or the symmetric group $S_{4}$. Since they are finitely generated, the invariant subgrassmannians in $\mathcal{G}_{k,n}(\mathbb{F}_{2}^{m})$ define a $t$-design where $t$ is the number of generators. The dihedral group $D_{m}$ is generated by a rotation $\tau_{m}$ of order $m$ and a reflection $\sigma$, $m=2\,s$, $s\geq 2$, then the corresponding invariant subgrassmannian in $\mathcal{G}_{k,n}(\mathbb{F}_{2}^{m})$ define a 2-$(m,s,\lambda)\,$design and so on for the other triangle groups. %where $\lambda$ is the number of 
%Moreover, for each conjugacy class $1, \sigma, \tau_{m}^{k},\tau_{m}^{s},\sigma\,\tau_{m}$, $k=1,\ldots,s-1$, in $D_{m}$ we get a 2-design, and so on for the other triangle groups.

%\cqd

\begin{lemma}\label{P1} The orbit of $\mathcal{G}_{2,4}(\mathbb{F}_{q})$ by the action of a rotation $\tau$ of angle $\alpha=\frac{2\pi}{n}$, corresponds to a cyclic code of order $n$. Furthermore, the orbit by the action of the element $\tau^{m}$ in $\GL(n,q)$ corresponding to the $m$-iterate composition of $\tau$ with itself is a quasi-cyclic code of index $\frac{m}{n}$.
\end{lemma}
{\it Proof.}
We study the action of a rotation element on the Grassmanian $\mathcal{G}_{2,4}(\mathbb{F}_{q})$ of lines in a 3-dimensional projective space $PG(3,q)$. We apply to any line $g$ % a symmetry $\sigma$ represented by the array of vectors $<(-1,0,0),(0,1,0),(0,0,-1)>$ and
 a rotation $\tau$ of angle $\alpha=\frac{2\pi}{n}$, represented by the array of vec\-tors $<(1,0,0), (0,cos(\alpha), sin(\alpha)), (0,-cos(\alpha), sin(\alpha))>$. It is easy to see that the orbit code by the composed action $ \tau^{m}$ with $m$  a divisor of $n$ is a quasi-cyclic code of index $\frac{m}{n}$.
 \cqd

%\subsection{Relation of $t$-designs with AG codes}
If we denote by $\mathcal{P}(n)$ the set of all linear subspaces inside the vector space $\mathbb{F}_{q}^{n}$, there is a natural metric on it defined by the function:
$$d_{S}(U,W):=\rm{dim}(U+W)-\rm{dim}(U\cap W).$$
The metric on $\mathcal{P}(n)$ induces a metric on the Grassmannian $\mathcal{G}_{n,k}(\mathbb{F}_{q})$. % With this metric it is possible to define a topology on the Grassmannian for which it is compact
For any subspace code $\mathcal{C}\subset \mathcal{P}(n)$, we define its distance through:
$$\rm{dist}\,(\mathcal{C}):=\rm{min}\,\{d_{S}(U,W)|\, U,W\in \mathcal{C},\, U\neq W\},$$ and its size as $M:=|\mathcal{C}|$.
A code is said to have minimum distance $d$ if $d_{S}(U,W)\geq d$ for all distinct words $U, W \in \mathcal{C}$. If the norm $|U|=w$ for every codeword in $\mathcal{C}$, then $\mathcal{C}$ is said to be of constant weight $w$. The number of codewords in $\mathcal{C}$ is called the size of the code.

We will say that a code $\mathcal{C}$ is of type $[n,k,d]$ if $\mathcal{C}$ has length $n$, minimum distance $d$, and its dimension is $k$. 

\begin{defi} Given a linear $[n,k,d]$-code, a parity check matrix for $\mathcal{C}$ is an $(n-k)\times n$ matrix $H$ of rank $n-k$ such that $\mathcal{C}=\{x\in (\mathbb{F}_{q})^{n}:\, Hc^{T}=0\}$. Then the dual code $\mathcal{C}^{\bot}$ is the linear $[n,n-k,d]$ code generated by the parity check matrix of $\mathcal{C}$. %Two linear codes of fixed parameters $[n,k,d]$ are called equivalent if they are permutation equivalent, that is if there is an element $g\in \GL(n,q)$ of rank $k$ and thus a linear transformation from $(\mathbb{F}_{q})^{n}$ to itself which permutes the underlying vector subspaces.
\end{defi}

\begin{defi} Given a linear $[n,k,d]$-code, a parity check matrix for $\mathcal{C}$ is an $(n-k)\times n$ matrix $H$ of rank $n-k$ such that $\mathcal{C}=\{x\in (\mathbb{F}_{q})^{n}:\, Hc^{T}=0\}$. Then the dual code $\mathcal{C}^{\bot}$ is the linear $[n,n-k,d]$ code generated by the parity check matrix of $\mathcal{C}$. %Two linear codes of fixed parameters $[n,k,d]$ are called equivalent if they are permutation equivalent, that is if there is an element $g\in \GL(n,q)$ of rank $k$ and thus a linear transformation from $(\mathbb{F}_{q})^{n}$ to itself which permutes the underlying vector subspaces.
\end{defi}
Any element $\sigma$ of the general linear group $\GL(n, \mathbb{F}_{q})$ induces another code $$\mathcal{C}^{\sigma} =\{\left( f(\sigma(\alpha^{i}))_{i=0}^{n}\right):\, f\in I\}.$$
A singer cycle of $\GL(n,\mathbb{F}_{q})$ is an element of order $q^{n}-1$. Singer cycles can be constructed, for example, by identifying vectors in $\mathbb{F}^{n}_{q}$ with elements of the finite field $\mathbb{F}_{q^{n}}$. Since multiplication by a primitive element $\alpha \in \mathbb{F}_{q^{n}}$ is a linear operation, it corresponds to a singer cycle in $\GL(n,\mathbb{F}_{q})$.
 For example, consider the codeword $a=(a_{1},\ldots, a_{n})\in \mathcal{C}$ and the permutation $\sigma \in \GL(n,\mathbb{F}_{q})$ which reverse the coordinates, so $\sigma(a)=(a_{n},a_{n-1},\ldots, a_{2},a_{1})$. If $\sigma(a)\in \mathcal{C}$, the code is called reversible.

The intersection code space $\mathcal{C}\bigcap \mathcal{C}^{\sigma}$ is given by the system of linear Diophantine equations:
 $$f(\alpha^{i})=g(\sigma(\alpha^{i}))\, \forall\, i=0,\ldots, n.$$
 In particular if $\sigma$ is a Singer cycle $\sigma(\alpha^{i})=(\sigma\alpha)^{i}$ permutes the elements $(\alpha^{i})$. Moreover if the permutation is in the automorphism group of the code, we get an equivalent code.
 Given two permutation codes $\mathcal{C}$ and $\mathcal{C}^{\sigma}$ the distance between them in the subspace metric is given by the formula:
$$d(\mathcal{C},\mathcal{C}^{\sigma}):= {\rm{dim}}(\mathcal{C}+\mathcal{C}^{\sigma})-{\rm{dim}}(\mathcal{C}\cap \mathcal{C}^{\sigma}).$$

\begin{defi}%We say that two $[n,k,d]$ binary codes $\mathcal{C}_{1}$ and $\mathcal{C}_{2}$ are equivalent if one can be obtained from the other one by permuting the $n$ symbols, that is, if there is a permutation $\pi$ such that $\pi(\mathcal{C}_{1})=\mathcal{C}_{2}$. 
Let $a$ be a word in $\mathcal{C}$, then the coset of $a$ is the set:
$$\{\pi(u)+ a: u\in \mathcal{C}, \pi\in S_{n}\}.$$
A coset leader is a word of minimum weight of any particular coset.

\end{defi}
% and adding a constant vector. Namely, there is a permutation $\pi$ and a vector $a$ such that $$\mathcal{C}_{2}=\{\pi(u)+ a: u\in \mathcal{C}_{1}\}.$$
%In particular, cosets of $C$ partition $\mathbb{F}_{q}^{n}$. 
Given a linear code $\mathcal{C}$, a non-trivial coset is a translation of $\mathcal{C}$ by a vector $v$ not in $\mathcal{C}$.
The main idea of coset coding is to map an information message not to a particular codeword but to a coset of this code.
We observe that two cosets are either equal or disjoint. %Actually to a syndrome associated with the coset
%Maximum security rate is known to be achieved when the code is a maximum distance code (MDC).
%\begin{example}

%\end{example}

\begin{remark} The incidence vectors of the blocks of a t-$(v,k,\lambda; q)$ design with maximum block intersection number $s$ form a constant weight code of weight $w=k$, lenght $n=v$, and minimum distance $d=2\,(k-s)$. %Genian Ge
\end{remark}

\begin{remark} A partition of the complete set of $k-$subspaces of $X$ into disjoint t-$(n,k,\lambda; q)$ designs is called a large set of $t$-designs over finite fields. Thus a partition of the Grassmannian $\mathcal{G}_{k,n}(\mathbb{F}_{q})$. Any point $\mathbb{F}_{q^{k},q^{n}}\subseteq \mathbb{F}_{q^{n}}$ in $\mathcal{G}_{k,n}(\mathbb{F}_{q})$ is a code of parameters $[n,k,d]$ where $d$ is the minimum distance defined as ${\rm{min}}\{d\,(V,W)|\, V,W \in \mathcal{G}_{q}(k,n),\, V\neq W\}.$
\end{remark}

\subsection{Relation of $t-$designs with AG codes}

Algebraic geometric codes (AGC), use as an alphabet a set $\mathcal{P}=\{P_{1}, \ldots, P_{N}\}$ of $N$-$\mathbb{F}_{q}$-rational points lying on a smooth projective curve $\mathcal{C}$ defined over  $\mathbb{F}_{q}$, that is, in projective coordinates $P_{i}=[a_{i}:1]$ with $a_{i}\in \mathbb{F}_{q}$. Namely, let $F/\mathbb{F}_{q}$ be the function field of the curve, $D$ a divisor of $F/\mathbb{F}_{q}$ supported on the set $\mathcal{P}$, and $G$ another divisor such that $\rm{Supp}\,G\cap \rm{Supp}\,D=\emptyset$. Then the geometric Goppa code $C(D,G)$ associated with the divisors $D$ and $G$ is defined by evaluation of a rational map $\varphi \in \mathcal{L}(G)$ in the linear series attached to the divisor $G$:
$$C(D,G)=\{(\varphi(P_{1}),\ldots, \varphi(P_{n})):\, \varphi \in \mathcal{L}(G)\}\subseteq \mathbb{F}_{q}^{n}.$$

It's an $\mathbb{F}_{q}-$subspace of $(\mathbb{F}_{q})^{n}$ and its dimension $k$ as an $\mathbb{F}_{q}-$vector space is the dimension of the associated Grassmannian code $G(n,k)$. Geometrically, it corresponds to a point in the Grassmannian $\mathcal{G}_{n,k}(\mathbb{F}_{q})$.

Observe that for the same subset of evaluation points and any $k'\leq k$, we have $G(n,k')\subseteq G(n,k)\subseteq \mathbb{F}^{n}_{q}$. In particular any $t-$design constructed from $G(n,k)$ is a $j-$design for $0\leq j \leq t-1$. It's well known that if ${\rm{deg}}(G)<n$, then $C(D,G)$ is a linear $[n,k,d]$ code over $\mathbb{F}_{q}$ with length $n$, $k=l(G)$ and minimum distance $d\geq n-{\rm{deg}}\,(G)$. We have $l(G)\geq deg\,(G)+1-g$ by the Riemann-Roch theorem, where $g$ is the genus of $F$.
% Referencia Niederreiter

More precisely, to each non constant rational function $\varphi$ over $C$ which is defined as the quotient of two polynomials $f(x), g(x) \in \mathbb{F}_{q}[x]$, one can associate a matrix $A$ with entries in the ring $\mathbb{F}_{q}[x]$. Then the generator matrix associated to the Goppa code $C(D,G)$ is defined to be the diagonal matrix with entries $q_{1}, q_{2},\ldots, q_{k}, k\leq n$, corresponding to the continued fraction expansion of the rational function $\varphi$. Namely, let us call $f_{0}:=f(x)$ and call $f_{1}$ the divisor polynomial $g(x)$, and $f_{2}$ the remainder polynomial, then by repeated use of the Euclid's algorithm, we construct a sequence of polynomials corresponding to the quotients $q_{1},\ldots, q_{k}$, $k\leq n$ of the continued fraction expansion $\frac{f}{g}=q_{1}+1/(q_{2}+1/q_{3}+1/(q_{4}+\ldots)))$.  Observe that if ${\rm{deg}}\,\frac{f(x)}{g(x)}<1$, then $q_{1}$ belongs to the ground field.These matrices are in correspondence with endormorphisms $T:\,R\rightarrow R$, of $\mathbb{F}_{q}[x]$-modules, where $R=\mathbb{F}_{q}(\alpha)$, and $\alpha$ is a generator of $\mathbb{F}_{q^{n}}$ as an $\mathbb{F}_{q}$-vector space. 
\begin{lemma}
The set of functions $\{q_{1},\ldots, q_{k}\}$ is in bijective correspondence with the set of codeword positions $f_{i}$ coming from the decomposition of the rational function $\varphi$ into partial fractions $f_{i}\in \mathbb{F}_{q}[x],\, 1\leq i \leq n$. Moreover, they  are linearly equivalently as $\mathbb{F}_{q}[x]$-vector spaces.\end{lemma}

{\it Proof.} We write the denominator $g$ of the rational function $\varphi$ as a product of powers of distinct irreducible polynomials. 
$$\frac{f}{g}=\frac{c_{1}}{x-\alpha_{1}}+\frac{c_{2}}{x-\alpha_{2}}+\ldots+\frac{c_{n}}{x-\alpha_{n}},$$ where the linear factors $(x-\alpha_{i})$ correspond to the roots of $g(x)$ counted with multiplicity. Since we are working over the finite field $\mathbb{F}_{q^{n}}$, the number of codeword positions $\{f_{i}\}_{i=1}^{n}$ is in correspondence with a base $\{1, \alpha, \ldots, \alpha^{n-1}\}$ of $\mathbb{F}_{q}^{n}$ as a vector space over $\mathbb{F}_{q}$, and they have the same cardinality as sets. \cqd

\begin{defi} Let $C$ be the AG code associated with a rational function $\varphi$ defined over a smooth projective curve $C$. A typical codeword is an element of the form $\sum_{j=1}^{k}a_{ji}f_{i} \equiv  0 \, \, {\rm{mod}}\,\, f$, where $f_{i}$ are the functions obtained from the decompostion of $f$ into partial simple fractions, and $a_{i}\in \mathbb{F}_{q^{n}}$.
\end{defi}

Let $\mathbb{F}_{l}$ be the finite field of prime-power order $l$ and let $F$ be an algebraic function field with full constant field $\mathbb{F}_{l}$. Note that $F$ can also be considered as an algebraic function field over any subfield $\mathbb{F}_{q}$ of $\mathbb{F}_{l}$.
One of the main problems in coding theory concerns the size of the alphabet $\mathcal{P}$, thus one of the aims is to obtain non trivial lower bounds of the number $N(F_{i})$ of rational places of towers of function fields $\{F_{i}/\mathbb{F}_{q}\}_{i=1}^{\infty}$ such that $F_{i}\subset F_{i+1}$.

The first case of study will be when the rational function $\varphi$ admits a decomposition into linear simple fractions. These rational functions define what are known as Reed-Solomon codes. In the case where $\alpha$ is a generator of $\mathbb{F}_{q^{n}}$, as an $\mathbb{F}_{q}$-vector space, the set of codeword positions is identified with the set of linear fractions 
$\{\frac{1}{x-1}, \frac{1}{x-\alpha}, \frac{1}{x-\alpha^{2}}\ldots, \frac{1}{x}\}$. Any linear combination of these elements produces a vector $a=(a_{1},\ldots, a_{n})$ in $(\mathbb{F}_{q})^{n}$ and thus a codeword of our AG code $\mathcal{C}$. Let $l$ be the maximum integer number such that the codeword position $a_{l}\neq 0$. The set of codewords $a$ satisfies the relation $$\sum_{i=1}^{n}a_{i}\frac{1}{x-\alpha_{i}}\cong 0\, {\rm{mod}}\, f,$$ where $L=\{\alpha_{i}\}_{i=1}^{n}$ is a subset of the Galois field $\mathbb{F}_{q^{n}}$. In several cases, it is possible to count the number of codewords of the AG code, by a simple count of the number of normalised polynomials of degree $l$ irreducible over $\mathbb{F}_{q^{n}}$. In the case of binary codes, where $q=2$, S. Bezzateev and N. Shekhunova (\cite{BS}) 
have obtained several closed formulas. The number of normalised polynomials $I_{2^{m}}(l)$ of degree $l$ over $\mathbb{F}_{2^{m}}$ satisfy the following equation:
\begin{equation}
I_{2^{m}}(l)=\frac{1}{l}\sum_{d / l}\mu(d)2^{m\, \frac{l}{d}},
\end{equation}
where $\mu(d)$ is the M\"oebius function. The number of unitary separable polynomials with coefficients from the field $\mathbb{F}_{2^{m}}$ whose degrees do not exceed $(l>1)$ is equal to:
\begin{equation}
N^{l}_{2^{m}}=\sum_{i=2}^{l}\, (2^{mi}-2^{m(i-1)})+2^{m}=2^{ml}.
\end{equation}

\begin{defi} The length of a codeword $(a_{1},\ldots, a_{n})$  in $(\mathbb{F}_{q})^{n}$ is $n=n_{1}+n_{2}+\ldots+n_{k}$, where $n_{i}$ is the number of positions of the vector $a$ with weight $v_{i}$ corresponding to the exponent of the corresponding fraction $f_{i}=\frac{1}{x-\alpha^{i}}$ in the partial fraction decomposition of the rational function $f$ associated to the AG code. 
\end{defi}

We observe that in the case of cyclic codes the weight $v_{i}$ coincides with the exponent of the corresponding function $f_{i}$ whose denominator is a linear function and thus with the integer $n_{i}$.

\begin{lemma}\label{equivalencia}
The set of AG codes defined over the Normal Rational Curve is in bijective correspondence with the set of generalised Reed-Solomon codes.% or cyclic codes.
\end{lemma}
{\it Proof.} We observe that the $n$-Veronese embedding of the $n$-dimensional projective space $\PG(n,q)$ maps the line spanned by the vector $v\in \mathbb{F}^{n+1}_{q}$ to the line spanned by $v^{n}\in \mathbb{P}S^{n}\mathbb{F}^{n+1}_{q},$ where $\mathbb{P}S^{n}\mathbb{F}^{n+1}_{q}$ is the projectivization of the $n-$tensor power of the vector space $\mathbb{F}^{n+1}_{q},$ which is a projective space of dimension $n$. In particular, if the finite field $\mathbb{F}_{q}^{n+1}$ is generated as a vector space over $\mathbb{F}_{q}$ by a unique element $\alpha\in \mathbb{F}_{q}$, then the set $\{1,\alpha, \ldots, \alpha^{n}\}$ forms a basis of $\mathbb{F}_{q}^{n}$. Thus the Normal Rational Curve is defined as:% where $v^{n}$
$$\mathcal{C}^{n}:=\{\mathbb{F}_{q}(1,\alpha, \ldots,\alpha^{n}):\, \alpha\in \mathbb{F}_{q}\cup \{\infty\}\}.$$
In other words, its underlying vector space is the $\mathbb{F}_{q}$-vector space whose elements are the polynomials of degree less than $n$ with coefficients in $\mathbb{F}_{q}$ that we will denote as $\mathbb{F}_{q}[x]_{n}$. Let $\alpha_{1},\ldots, \alpha_{n}$ be a sequence of $n$ distinct elements in $\mathbb{F}_{q}$, if $k\leq n$, then the map
\begin{equation}
\epsilon: \mathbb{F}_{q}[x]\rightarrow \mathbb{F}_{q}^{n}, \ \ f\mapsto (f(\alpha_{1},\ldots, \alpha_{n}))
\end{equation}
is injective, since the existence of a non-zero polynomial of degree less than $k$ vanishing on all $\alpha_{i}$ implies $n<k$ by the fundamental theorem of algebra (a non-zero polynomial of degree $r$ with coefficients in a field can have at most $r$ roots). The image of $\epsilon$ is therefore an AG code of type $[n,k,d]$, where the minimum distance $d$ is always at least $n+2$. Just observe that since $\mathcal{C}^{n}$ is a  Normal Rational Curve in $\mathbb{P}({\mathbb{F}_{q}^{n}})$, any $n+1$ of its points happen to be in general position.

Reciprocally, the AG codes of dimension $n$ defined over $\mathcal{C}^{n}$ are constructed by evaluating non-zero polynomials of degree less than $n$ over a location set $\{1,\alpha_{1}, \ldots, \alpha_{n}\}$ which coincide with the set of Reed-Solomon codes. Namely, consider a Reed-Solomon code of parameters $[n,k,d]$ over a finite field $\mathbb{F}_{q}$, with parity check polynomial $h(x)=\prod_{i=1}^{q}(x-\alpha^{i})$, where $\alpha$ is a primitive root of $\mathbb{F}_{q}$ such that $\alpha^{k+1}=\alpha+1$. Any codeword $(c_{0},c_{1},\ldots, c_{n-1})$ can be expanded into a $q$-ary $k$ vector with respect to the basis $\{1, \alpha,\ldots, \alpha^{k-1}\}$, that is, codewords from an $[n,k,d]$ code  over a finite field are identified with the coefficients of a degree $k-1$ polynomial $f(x)\in \mathbb{F}_{q}[x]$.
\cqd

%As a consequence of Lemma \ref{equivalencia}, codewords from an $[n,k,d]$ over a finite field are identified with the coefficients of a degree $k-1$ polynomial $f(x)$.% under the isomorphism $(a_0},\ldots, a__{n-1})\mapto a_{0}+a_{1}x+\ldots+a_{n-1}x^{n-1}.$

%Let $F$ be a field of characteristic $p$, in particular it contains a copy of $\mathbb{F}_{p}$. Given integers $n, k$ and $d$ with $1\leq k\leq n$, an $[n,k,d]_{F}-$code over $F$ is a subspace $C$ of $F^{n}$ of dimension $k$, such that every nonzero codeword $0\neq \alpha \in C$ satisfies $wt(\alpha)\geq d$, where the weight of a vector $\alpha=(a_{0},\ldots ,a_{n-1})\in F^{n}$ is the number of non-zero components $a_{i}$. 

\begin{example} Consider the AG code defined by the rational function $G(x)=\frac{5x^{2}+20x+6}{x^{3}+2x^{2}+x}$ which admits as decomposition into partial fractions %the expression 
$G(x):=\frac{6}{x}-\frac{1}{x+1}+\frac{9}{(x+1)^{2}}$. The presence of a double factor $(x+1),,,,,^{2}$ corresponds to the existence of an eigenspace $E$ in the vector space $\mathbb{F}_{q}^{n}$ and thus an $\alpha-$splitting subspace where the operator $\alpha$ is just the linear operator $A-\lambda I$, with $\lambda$ the eigenvalue of multiplicity 2  associated to $E$ and $A$ is the generator matrix of the code.
\end{example}

\begin{prop}\label{fields}
%Let $F$ be a field of characteristic $p$. 
The variety of $[n,k,d]_{q}$-codes over $\mathbb{F}_{q}$ is parametrized by a Grassmannian $\mathcal{G}_{n,k}(\mathbb{F}_{q})$ of $k$-di\-men\-sional subspaces in the $\mathbb{F}_{q}$-vector space $\mathbb{F}_{q}^{n}$, and the set of Reed-Solomon $(RS)$ codes arises as the set of $S_{n}$-invariants.% so in particular, as the Grassmannian is a compact variety, the set of RS codes is a closed set in the Zariski topology.
\end{prop}
{\it Proof.}
Let $n=r\,s$ be a factorisation of an integer positive number $n$ into irreducible coprime factors and assume $s<r$, then there is a sequence of field extensions $\mathbb{F}_{q^{r}}\subset \mathbb{F}_{q^{s}}\subset \mathbb{F}_{q^{n}}$.
Namely, consider the map
 $T_{n}:\,   F^{n}  \mapsto   F^{n}$
$$t_{j}=(-1)^{j}\sigma_{j}(x_{1},\ldots, x_{n}),$$ where $\sigma_{j}$ is the $j^{th}$ elementary symmetric function in the variables $x_{i}$. Thus $\{t_{j}, j=1,\cdots n\}$, are the coefficients of the equation:
$$f(z,t_{1},\ldots, t_{n})=z^{n}+(-1)\,t_{1}z^{n-1}+\cdots +(-1)^{n}\,t_{n}=$$ $$(z-x_{1})\,(z-x_{2})\cdots (z-x_{n}).$$ Then by Hilbert's irreducibility theorem (see Theorem 1 of \cite{ Se}), it is well known that the splitting field of the polynomial $f(x)=x^{n}-t_{1}x^{n-1}+\ldots+(-1)^{n}t_{n},$ is the field of $S_{n}$ invariants of the polynomial $f(z,t_{1},\ldots, t_{n})$, where $S_{n}$ is the symmetric group  in $n$ variables and %\in F[z,x_{1},\ldots,x_{n}]$ 
it contains an extension $\mathbb{F}_{q^{n}}$ of $\mathbb{F}_{q}$.  Moreover, for any divisor $r$ of $n$, one can consider the field of $S_{r}$ invariants, and apply Hilbert theorem to the symbols $\alpha, \alpha^{q^{2s}},\ldots, \alpha^{q^{rs}}$, where $n=rs$. Then we get an extension $\mathbb{F}_{q^{s}}$ of $\mathbb{F}_{q^{r}}$ and all its $\mathbb{F}_{q}$-subspaces are stable under $Gal(\mathbb{F}_{q^{s}}/\mathbb{F}_{q^{r}})$. These are just the RS codes.%These are known as generalised Reed-Solomon codes.

\cqd
\begin{coro} The set of RS codes is a closed set in the Zariski topology.
\end{coro}
{\it Proof.}This follows easily as a consequence of Proposition \ref{fields}, since the Grassmannian is a compact variety. It is well known that the corresponding points 
$\mathbb{F}_{q^{k},q^{n}}\subset \mathbb{F}_{q^{n}}$ and $\mathbb{F}_{q^{n-k},q^{n}}\subset \mathbb{F}_{q^{n}}$ in the Grassmannians $\mathcal{G}_{k,n}(\mathbb{F}_{q})$ of $k-$dimensional subspaces and the Grassmannian $\mathcal{G}_{n-k}(\mathbb{F}_{q})$ of $n-k$ dimensional subspaces are respectively dual subspaces in the underlying vector space $(\mathbb{F}_{q})^{n}$ for the Euclidean inner product. Note that the Hamming weight is preserved under invertible linear transformation
\cqd

\begin{theorem}\label{goppa}
Let  $S\leqslant \GL(n,\mathbb{F}_{q})$ be a subgroup containing a primitive element $\alpha$ in $\mathbb{F}_{q}$, where $q\geq 2$.
%Let $\alpha\in \mathbb{F}_{q^{n}}$ be a generator of the underlying vector space over $\mathbb{F}_{q} $. 
%Whenever $\alpha$ is an element of an invariant subgroup  $S\leqslant \GL(n,\mathbb{F}_{q})$, 
Then the family of AG codes $\{C^{\sigma}\}_{\sigma \in S}$ constitute a j-$(v,r,\lambda)$ design where $j$ is the number of generators of the subgroup $S$, $r$ is the number of orbits in the Grassmannian $\mathcal{G}_{k,n}(\mathbb{F}_{q})$ by the action of the subgroup $S$, $v$ is the size of the code and $\lambda$ is the number of $\alpha-$splitting subspaces of $\mathbb{F}_{q^{n}}$.

\end{theorem}
{\it Proof.} 
Consider the family of AG codes constructed out of the vector space of polynomials $I=\{f\in \mathbb{F}_{q}[x]: \partial\, f\leq k\}$, where $\partial\, f$ is the degree of the polynomial and fix a basis $\{1,\alpha, \alpha^{2},\ldots, \alpha^{n-1}\}$. To each polynomial $f$ we associate the AG code $C=
 \{\left( f(\alpha^{i})_{i=0}^{n}\right):\, f\in I\}$. 
 
We construct a $j$-design where the point set are the codewords of the AG codes and each AG code of constant dimension $r$ is a block, where $r$ is the number of orbits in the Grassmanian $\mathcal{G}_{k,n}(\mathbb{F}_{q})$ by the action of any representative of a conjugacy class in $S$. %which coincides with the number of conjugacy classes in S.
% the general linear group $\GL(n,\mathbb{F}_{q})$ considered in (\ref{action}). 
Each polynomial $f(x)=a_{0}+a_{1}x+\ldots+a_{j}x^{j}$ defines a codeword $(a_{0},\ldots,a_{j})\in \mathbb{F}_{q}^{j}$ of the code. 
Since the number of invariant polynomials  in $I$ by  conjugated elements $A$ and $B$ in $\GL(n,\mathbb{F}_{q})$ is the same, $r$ is the number of conjugacy classes in $S$. The intersection vector space is given by the evaluation set 
$$\{\left(f(\alpha^{i})\right)_{i=0}^{n}:\, f\in \mathbb{F}_{q}[x]_{j},\, 1\leq j\leq r\},$$
of the polynomials of degree $j$ with coefficients in $\mathbb{F}_{q}$, $v$ is the size of the block codes which is constant and $\lambda$ is the number of $\alpha$-splitting subspaces of $\mathbb{F}_{q^{n}}$ as computed in \cite{BM2}.
 \cqd

Let $F$ be a field of characteristic $p$ and $\alpha \in \overline{F}$ be an $n^{th}$ primitive root of unity, where $\overline{F}$ denotes the algebraic closure of $F$. %$m_{i}(x)$ be 
The $n^{th}$ cyclotomic polynomial $\Phi_{n}(x)=\prod_{1<j<n, (j,n)=1}\,(x-\alpha^{j})\in \overline{F}[x]$ is the minimal polynomial of $\alpha$ over $F$. It is monic of degree the Euler's totient function $\varphi(n)$.
It has integer coefficients and it is irreducible over $\mathbb{Q}$. In $\mathbb{Q}[x]$, we have the factorization into irreducible polynomials:
$$x^{n}-1=\prod_{d|n}\Phi_{d}(x).$$
By M\"oebius inversion:
$$\Phi_{n}(x)=\prod_{d|n}(x^{d}-1)^{\mu(n/d)}$$ 
\begin{example} Consider the action of a permutation matrix $\beta$ in $\GL(n,q)$ given by the element $\beta\in \mathbb{F}_{q}^{n}:\,x\rightarrow a\,x$, which is given by multiplication by an element $a\in \mathbb{F}^{*}_{q}$ of multiplicative order $k>1$ with $n=km$. Then by Lemma \ref{equivalencia}, counting the codewords of the constant dimension code defined by the action of $\beta$ in the Grassmannian $\mathcal{G}_{k,n}(\mathbb{F}_{q})$, is equivalent to count the number $\mathcal{N}_{a,m}$ of irreducible monic polynomials of degree $n$ such that $f(x)=f(ax)$. This number is expressed in terms of the  M\"oebius function:
$$\mathcal{N}_{a,m}=\frac{\Phi(k)}{km}\sum_{d|m, gcd(d,k)=1}=\mu(d)(q^{\frac{m}{d}} -1).$$
\end{example}
\begin{example}
We consider the roots of the polynomial $x^{8}-1\in \mathbb{F}_{5}[x]$ in the splitting field $\mathbb{F}_{5^{2}}$. The decomposition into irreducible polynomials over $ \mathbb{F}_{5}[x]$ is $(x-1)(x+1)(x-2)(x+2)(x^{2}+1)(x^{2}-2)(x^{2}+2)$. Now, we consider the field extensions $F_{1}:=\mathbb{F}_{5}[x]/(x^{2}-2)$ and $F_{2}:= \mathbb{F}_{5}[x]/(x^{2}+2)$ of $\mathbb{F}_{5}$ that are isomorphic to the field extension $\mathbb{F}_{25}$ of $\mathbb{F}_{5}$. Call $\alpha$ the root of $x^{2}-2$ in the field extension $F_{1}$, then $4\cdot \alpha$ is the other root of $x^{2}-2$, and $2\cdot \alpha$, $3\cdot \alpha$ the roots of $x^{2}+2$ in $F_{1}$. So $g(x)=(x-\alpha)(x-2\alpha)(x-3\alpha)(x-4\alpha)$ generates a Reed-Solomon code %or cyclic code 
over $\mathbb{F}_{5}[x]/(x^{8}-1)$.
We say that two roots are conjugated if they are roots of the same polynomial in the decomposition of $x^{8}-1$ in $\mathbb{F}_{5}[x]$, in particular this defines a non-crossing partition of the total set of roots and it is a 2-design of the splitting field $\mathbb{F}_{25}$ with parameters $n=2$, $k=4$ and $\lambda=2$.

From a geometric point of view, %if we consider the factorisation of the polynomial $(x^{2}-1)(x^{2}+1)(x^{4}+1)$ over $\mathbb{F}_{5}[x]$, we see that 
the point $(\alpha, 0)\in \mathbb{P}(\mathbb{F}_{q}^{2})$ with $\alpha^{4}=4$ is an $\mathbb{F}_{25}-$rational point of the affine curve $y^{2}=(x^{4}+1)$. The other rational places are
$(2,0), (-2,0)$ and the place $(0,\alpha)$ at $\infty$. 
\end{example}

It is well known how to factorize a polynomial over finite fields (see for example \cite{PFG}). In \cite{BM2} we give an updated proof expressing the number of polynomials decomposable into distinct linear factors in terms of Stirling numbers.

%\begin{theorem}\label{count} The number of  polynomials decomposable into distinct linear factors over a finite field $\mathbb{F}_{q^{n}}$ of arbitrary characteristic  a prime number $p$, is equal to
%$\sum_{k=1}^{n}(q)_{k}$, where $(q)_{k}$ is the falling factorial polynomial 

%\noindent $q\cdot (q-1)\ldots (q-(k-1))$, divided by the order of the affine transformation group of $\mathbb{A}^{1}=\mathbb{P}^{1}\backslash \infty$, that is $q^{2}-q$.

%\end{theorem}

%{\it Proof.} We need to count all the polynomials $f_{n}(x)$ in one variable of degree $n$ fixed. We assume that our polynomial $f_{n}(x)$ decomposes into linear factors, otherwise we work over $\bar{\mathbb{F}}_{q}[x]$, where $\overline{\mathbb{F}}_{q}$ denotes the algebraic closure of the finite field $\mathbb{F}_{q}$. Since the number of ordered sequences on $q$ symbols is $q!$, it follows that the number of monic polynomials with $n-1$ different roots is $q(q-1)(q-2)\ldots (q-n+1):=(q-2)_{n}=\sum_{k=0}^{n}s(n,k)\,q^{n}$, where $s(n,k)$ is the Stirling number of the first kind and it counts the number of ways to partition a set of cardinality $n$ into exactly $k$ non-empty subsets. Now we observe that polynomials are invariant by the action of automorphisms of the affine line, so we must divide this number by the order of this group which is $q^{2}-q$.
%\cqd

%As an application %of Theorem \ref{count}, 
Given an integer $n$, it is possible to count the number of cyclic codes of parameters $[n,k]$ for each $0\leq k \leq n$ and set of roots $\alpha_{1},\ldots, \alpha_{k}$ in the splitting field of $x^{n}-1$, the corresponding polynomial $g(x)=\prod_{i=1}^{k}(x-\alpha_{i})$ generates a linear cyclic code in the ring $\mathbb{F}_{q}[x]/(x^{n}-1)$. Thus for each $0\leq k\leq n$ there are exactly $(q)_{k}/(q^{2}-q)$ cyclic codes.
These codes are of great importance in ADN-computing and as they are linear codes, they can be described as function fields. %Let $\alpha$ be a primitive element of the underlying vector space over $\mathbb{F}_{q}$. Since the base field is of characteristic $p$, $x^{n}-1$ has $n$ different zeroes. Let $\bar{\mathbb{F}}_{q}[x]$ be the extension field containing the $n^{th}$ roots of unity $1,\alpha, \ldots, \alpha^{n-1},$ where $\alpha^{n-1}+\alpha^{n-2}+\ldots+\alpha+1=0$. Moreover the set $\{1,\alpha, \ldots, \alpha^{n-1}\}$ constitutes a basis over the prime field $\mathbb{F}_{p}$, and the field extensions $\mathbb{F}_{p^{n}}\cong \mathbb{F}_{p}[x]/(x^{n-1}+\dots +x+1)$ are isomorphic.

\begin{remark}
A much greater variety of linear codes is obtained if one uses places of arbitrary degree rather than just places of degree 1 as in Goppa's construction. For example, the polynomial $x^{3}+4$ factorises as $(x-1)(x^{2}+x+1)$ over $\mathbb{F}_{5}[x]$, then the roots of the polynomial in the splitting field $\mathbb{F}_{5}[x]/x^{2}+x+1 \cong \mathbb{F}_{25}$ correspond to one place of degree 2 over the function field $\mathbb{F}_{5}(x)$ but of degree 1 over $\mathbb{F}_{25}$.

\end{remark}

\subsection{$t-$designs with an action of a $p-$group.} %over $\mathbb{F}_{5}[x]/(x^{8}-1)$.}
During the last years, there has been an increasing interest in studying finite abelian groups due to its relationship with public key cryptography, quantum computing and error-correcting codes. Abelian groups as the groups $\mathbb{Z}_{n}^{*}$ of invertible elements of $\mathbb{Z}_{n}$, multiplicative groups of finite fields, the groups of elements of elliptic curves over finite fields, finite $p-$groups with unique cyclic subgroups of given order have been used for the designation of public key cryptosystems. In order to use cryptography to insure privacy, it is currently necessary for the communicating parties to share a key which is known to no one else. %A private conversation between two people is a very common occurrence in business, however, it is unrealistic to expect initial business contacts to be postponed long enough for keys to be transmitted by some physical means.%referencia
As we showed in Theorem \ref{goppa}, we can construct $t$-designs from any $p-$group containing a cyclic subgroup. 

\begin{prop}
For $q\geq 2$, the group $\GL(n,q)$ contains a least two different cyclic subgroups of orders $q-1$ and $q+1$ respectively. Each one corresponding to elements $\alpha, \gamma$ in $\GL(n,q)$ fixing  an $\alpha-$splitting and $\gamma-$splitting subspaces respectively and the $t-$designs whose incidence vectors are the $\alpha-$splitting subspaces and $\gamma-$splitting subspaces respectively correspond to RS codes of length $n=q-1$ (respectively $n=q+1$) and dimension $r$ the maximum divisor of $n$.
\end{prop}
{\it Proof.} %By proposition \ref{lambda}, f
For any divisor $d$ of $q-1$ (respectively of $q+1$), the $d$-$(q-1,\frac{q-1}{d},\lambda)$ design (respectively $d$-$(q+1,\frac{q+1}{d}),\lambda)$) corresponds to RS codes of length $n=q-1$ (respectively $n=q+1$) and dimension $d$. Moreover the matrix  $A$ of row vectors the incidence vectors of the design, satisfies $\frac{q+1}{r}\leq rank(A) \leq \frac{q-1}{2}$, where $r$ is the maximum divisor of $q-1$ (respectively $q+1$).  \cqd

\begin{remark} 
The generators of the cyclic groups of order $q-1$ ($q+1$, respectively) are the relative integers coprime with $(q-1)$ (respectively with $q+1$),  that is $\varphi(q-1)$ (respectively $\varphi(q+1)$).
 Let $m=\varphi(q-1)$ (respectively $m=\varphi(q+1)$),  by Theorem \ref{goppa}, the family of RS codes of length $q-1$ (respectively $q+1$) constitute a $m$-$(q-1,r,\lambda)$ design (respectively a $m$-$(q+1,r,\lambda)$ design). These codes are indeed AG codes arising from genus 0 curves, and by Riemann-Roch theorem, their parameters satisfy the bound $d\geq n+1-k$, where $d$ is the minimum distance. % in the Hamming metric.

\end{remark}

%\vspace{0.5cm}
The normalizer groups of the cyclic groups generated by $\alpha$ and $\gamma$ are dihedral groups, and it is possible to construct $t$-designs from them as we showed in Proposition 3.5 of \cite{BM2}.%\ref{PropSteiner}. %point $a$ and $b$ 
Next tables show $t$-designs constructed from abelian $p$-groups and their normalizers. %Some of these designs were constructed from linear groups $GL(2,q)$, $q\leq 23$. we construct transitive t-designs from the linear groups L(2, q), q ? 23.
%\'ic
%D. Crnkov\'ic and A. Svob construct in \cite{CS} some of these designs from linear groups with $q\leq 23$. 
We assume that $q\leq 31$.
%\vspace{1cm}
\newpage
\begin{center}
Table 1. $t-$designs constructed from a $p$-group.
\end{center}
\begin{center}
%\hline Table 1
\begin{tabular}{|c|r|r|r|}
%\hline Table 1 & & \\
\hline
   $q$ & group type & $t-(n,k,\lambda)$   \\
\hline 
  $q$ odd prime $q+1\equiv 0 (3)$ & cyclic of order $q+1$ & 3-$(q+1, \frac{q+1}{3},\lambda)$ \\
  
   $q$ odd prime $q-1\equiv 0(3) $ & cyclic of order $q-1$ & 3-$(q-1, \frac{q-1}{3},\lambda)$ \\

$q$ odd & cyclic of order $q$ & 3-$(q,\frac{q}{3},\lambda)$ \\

 $q=p^{e}$ & abelian $p$ group & $p$-$(q, p^{l},\lambda)$, $l<e$  \\
  \hline    
 \end{tabular}
 \end{center}
 
 %\newpage
 \vspace{0.5 cm}
\begin{center}
Table 2. $t-$designs constructed from their corresponding normalizers.
\end{center}
 \begin{center}
\begin{tabular}{|c|r|r|}
\hline
   $q$ & group type & $t-(n,k,\lambda)$   \\
\hline  $q$ odd & dihedral of order $2\,(q-1)$ & 3-$(2\,(q-1), q-1,\lambda)$ \\

   $q$ even & dihedral of order $2\,(q+1)$ & 2-$(2\,(q+1), q+1,\lambda)$ \\

 $q=p^{e}$ & Borel  & $p$-$(q\,(q-1),(q-1),\lambda)$ \\
 \hline
 \end{tabular}
 \end{center}
\vspace{0.5cm}

Recall that the dihedral group of order $2\,(q-1)$ (respectively of order $2\,(q+1)$) is generated by a rotation $\tau_{q-1}$ of order $(q-1)$, (respectively of order $2\,(q+1)$ and a reflection. In particular the discrete logarithm problem (DLP) applied to this group reads:
Given an element $h\in D_{2\,(q-1)}$ find an integer $m$ satisfying $\tau^{m}=h$. The smallest integer $m$ satisfying the identity is called the index of $h$ with respect to $\tau$, and is denoted as $m=log_{\tau}(h)$ or $m=ind_{\tau}(h)$.
The DLP is used as underlying hard problem in many cryptographic constructions, including for example Diffie-Hellman key exchange, \cite{DH}. %and digital signatures.
%\begin{example}
Solving DLP takes time that is exponential in the order of the group $G$. For example, the group defined by the elliptic curve over a finite field $\mathbb{F}_{p}$ takes time $O(\sqrt(p)$. For this reason, it is used for cryptographic purposes. 

\subsection{Diffie-Hellman key exchange for dihedral groups}
%A public-key encryption system is a triple of algorithms $(G,D,E)$, the generating algorithm which produces the public key $ pk$ and the secret key $sk$. The encryption algorithm takes the public key and a message and produces a cipher tex. The cipher tex is fed into the decryption algorithm and using the secret key it produces the corresponding message $m$.
%It is consistent, that is, $\forall (pk, sk)$ output by $G$:
 %$\forall m \in M$,  $ D(sk, E(pk,m))=m$. % that is, 
%Basically it says that 

% If we encrypt a message with a given public key and then we decrypt under the secret key we obtain the original message.
% The public key  defines a bijective function $F(pk,.):\, X\rightarrow X$ from the space of messages $X$ to itself and the secret key  defines the inverse function $F^{-1}(sk,.)$.
 %%from the set $Y$ to the set $X$ that inverts $F(pk,.)$. 
 %More precisely, if we look at a key pair $(pk, sk)$ generated by the generating algorithm $G$, then it happens that if we evaluate the function at the point $x$ and then we apply the inverse function, we get the original message back, that is:
% $\forall (pk,sk)$ output by $G$, $\forall x\in X:  \ \  F^{-1}(sk, F(pk,x))=x$. Then $F$ is called a trapdoor function.

In Diffie-Hellman key exchange cryptosystem, the public key is an element of a group $G$ of public knowledge, in the case of study is a dihedral group of order $2\,(q-1)$ generated by a reflection $\sigma$ of order 2, and a rotation $\tau$ of order $(q-1)$. The generating algorithm produces an element $\tau$ which is the public key. Observe that since the group generated by $\tau$ is cyclic of order $q-1$, its elements $1,\tau, \tau^{2},\ldots, \tau^{q-2}$ are roots of unity, that is, $x-\tau^{i}$ divides the polynomial $x^{n}-1$, and the code is cyclic. Moreover it is the RS code of length $n=q-1$ and dimension $k$. %This is an example of RS code-based public encryption scheme.
There are two participants involved in the encryption process, participant $P_{1}$ and participant $P_{2}$. %Each one sends messages to the other enciphered in the receiver's public enciphering key and enciphers the messages received using his own secret deciphering. %
A third party eavesdropping on this exchange must find it computationally infeasible to compute the key from the information overheard. %Each user generates a pair of inver%Diego and Elisa and the trapdoor function is the discrete logarithm.
First par\-ti\-ci\-pant $P_{1}$ randomly choose a secret $0<d<(q-1)$ and computes $D=\tau^{d}$. Second participant $P_{2}$ randomly choose a secret $0<e<(q-1)$ and computes $E=\tau^{e}$. Participant $P_{1}$ sends $D$ to participant $P_{2}$, and $P_{2}$ sends $E$ to $P_{1}$. Then $P_{1}$ computes $E^{d}=\tau^{e\,d}$ and $P_{2}$ computes $D^{e}=\tau^{d\,e}$ so that both participants $P_{1}$ and $P_{2}$ have the shared value $\tau^{d\,e}$ up to reflection $\sigma\in D_{q-1}$. Thus computing $\tau^{e\,d}$ from $\tau^{e}$ requires solving the DLP $log_{\tau^{e}}(\tau^{e\,d})=d$. % that is, since the construction is symmetric with respect to the variables $a$ and $b$, 
%there is a reflection $\sigma\in D_{q-1}$ such that $\sigma\, (\tau^{ab})=\tau^{ba}  D_{q-1}$.
%\end{example}
At each time $t$, the probability of select the element $\tau^{j}$ is distributed as a Bernoulli distribution with $Pr(x=\tau^{j})=\left(\frac{1}{q}\right)^{j}(1-\frac{1}{q})^{j}, \, j\in \{0,1,\ldots q-1\}.$
%There are two people involved in the encryption process, Bob and Alice and the trapdoor function is the discrete logarithm. The cryptographic protocol consist in repeating $k$ times the following 3 exchanges:

%\begin{center}
%\begin{tabular}{|c|}
%\hline
%Bob choose a secret $0<d<(q-1)$ and computes $D=\tau^{d}$.  \\
% Alice choose a secret $0<a<(q-1)$ and computes $A=\tau^{a}$.  \\

%Bob sends $D$ to Alice, and Alice sends $A$ to Bob.  \\

%Then Bob computes $A^{d}=\tau^{a\,d}$ and Alice computes $D^{a}=\tau^{d\,a}$ so that   \\

%Bob and Alice have the shared value $\tau^{d\,a}$ up to reflection $\sigma\in D_{q-1}$.  \\
%\hline

% \end{tabular}

 %\end{center}

%The network model
\section{Relation of $t$-designs with graph networks} 
%Methods of classification of data sets which can be recorded as strings of sequences of letters over a finite alphabet, use graph theory. % (see \cite{KRY}). 
A standard way to do natural language processing (NLP) are networks. 
   A sample graph is a network constructed from data that has been collected from a random sample of nodes. Many problems are motivated by the need to infer global properties of the parent network (population) from the sampled version. Counting the number of features in a graph is an important statistical and computational problem. These features are typically basic local structures like motifs or graphlets (e.g., patterns of small subgraphs). 
One of the most important classes of graphs considered in this framework is that of Cayley graphs. 
Consider a network represented by a directed multigraph
%Consider a network graph 
$G=(V(G),E(G))$, with vertex set $V(G)$ and edge set, $E(G)$ with error free unit capacity edges, that is, a graph with loops (edges whose endpoints are equal) and multiple edges. % where two vertices are set to be adjacent if there is an edge between them.
%There are several source nodes and several destination nodes and data is transferred over a network using packets. % A packet is a $m$-length vector over a finite field $\mathbb{F}_{q}$.
 A subset of the vertex set is called independent set, if there is no edge between vertices in $X$. A matching is a set of disjoint edges of a graph.
A clique in an undirected graph is a subset of its vertices such that every two vertices in the subset are connected by an edge. Let $G$ be a simple graph and $H$ a subgraph of $G$. A $G$-design of $H$ is a pair where $X$ is the vertex set of $H$ and $\mathcal{B}$ is an edge-disjoint decomposition of $H$ also known as partition of the vertex set. %Each block of a given partition identifies the positions where the $L-$matchings take place.
We say $G$ is a split graph if the vertex set $V(G)$ can be partitioned into a clique $C$ and an independent set $I$, where $(C,I)$ is called a plot partition of $G$.
The best known cryptographic problem is that of privacy: preventing the an authorised extraction of information from communications over an insecure channel. In order to use cryptography to insure privacy, however, it is currently necessary for the communicating parties to share a key which is known to no one else.
For applications in data privacy, we are interested to identify subgraphs in the complete graph that are identical, this gives a measure of the degree of anonymity of the graph. 

An automorphism is a permutation of the vertices of the network which preserves adjacency. The set of automorphisms under composition forms a group $\rm{Aut}(G )$ of size $a_{G}$ which compactly describes network symmetry. The orbit of a vertex $v\in V(G)$ is the set:
$$\triangle(v)=\{\pi v\in V(G):\, \pi\in \rm{Aut}(G)\}.$$
Automorphism group orbits naturally partition network vertices into disjoint structural equivalence classes. Since two vertices in the same orbit may be permuted without altering network adjacency, they are structurally equivalent in the strongest possible way: they play exactly the same structural role in the network. %The network redundancy is defined by the quantity $r_{G}=\frac{N_{Q}-1}{N_{G}}$, where $N_{\mathcal{Q}}$ is the number of network orbits and $N_{G}$ is the number of vertices in the network, \cite{KTSZ}.%referencia

Let $G$ be a network with automorphism group $\rm{Aut}(G)$. Let $1\neq S$ be a set of generators of $\rm{Aut}(G)$. Suppose that we partition $S$ into $n$ support-disjoint subsets $S=S_{1}\cup \ldots \cup S_{n}$ such that each $S_{i}$ cannot itself be decomposed into smaller support-disjoint subsets. Call $H_{i}$ the subgroup generated by $S_{i}$. Since $S$ is a generating set and elements from different factors $H_{i}, H_{j}$ commute, this procedure gives a direct product decomposition: $$\rm{Aut}(G)=H_{1}\times H_{2}\times\ldots \times H_{n}.$$
The network automorphism group decomposition relates automorphism group structure to network topology. Moreover, automorphism groups of real-world networks such as scientific collaboration networks or technological networks such as the internet can typically be decomposed into direct and wreath products of symmetric groups, (see \cite{BRJ}).
Reciprocally, given a group presentation $S$ of $\GL(n,\mathbb{F}_{q})$, we can attach to it a Cayley graph which is defined as the directed graph having one vertex associated with each group element and directed edges  $(e_{1}, e_{2})$ whenever $e_{1}e_{2}^{-1} \in S$. The Cayley graph may depend on the choice of a generating set, and it is connected if and only if $S$  generates $\GL(n,\mathbb{F}_{q})$.

A subset $S$ of an additive group is called sum-free if it contains no elements $x, y, z$ such that $x+y=z$. In particular, this means the corresponding vertices constitute an independent set in the Cayley graph, (see \cite{KTSZ}). Moreover, there are two distinguished sets of vertices, the set of independent vertices $L$ which satisfy the property that there is no edge between two vertices and the complement graph.%graph 

We pass from network topology to vector network coding by thinking of the vertex set $V(G)=\{a_{1},\ldots,a_{n}\}$ as an alphabet of $n$ letters and defining a vector space $V$ on these $n$-generators over a ground field $k$. There is a natural representation $\rho: S_{n}\rightarrow \GL(V)$, where $S_{n}$ is the group of permutations of $n$ elements. As we showed in \cite{BM1}, vector network coding and moreover codes over a finite field $\mathbb{F}_{q}$ are very much related with the study of the representation theory of the symmetric group over finite fields and further with the representation theory of $\GL(n,\mathbb{F}_{q})$ over finite fields. From this representation, one can recognise more easily patterns and extract information from them.
%As we showed in \cite{BM1} to study Grassmannian codes over a finite field $\mathbb{F}_{q}$ is reduced to study the representation theory of $\GL(n,\mathbb{F}_{q})$. This is clear, as we are interested to find the matrix codifying that it is the matrix of an endomorphism of $\mathbb{F}_{q}[x]-$modules.  From this re\-pre\-sen\-ta\-tion, one can recognise more easily patterns and extract information from them. 
In terms of designs over $\mathbb{F}_{q}$, we want to understand which subspaces are invariant by the action of elements of the general linear group $\GL(n, \mathbb{F}_{q})$ or finite subgroups of $\GL(n,\mathbb{F}_{q})$. In this way, one can construct designs with prescribed groups where the blocks are the orbits by the action, and thus to generalise to other Galois extensions not necessarily cyclic. 

%For other general linear groups, we study generalised Grassmannians or more commonly known as flag varieties.
%Fix a partition $\lambda=(\lambda_{1},\ldots, \lambda_{r})$ of $n$ and let $\mathcal{F}_{\lambda}=\mathcal{F}_{\lambda}(\mathbb{F}_{q})$ be the variety of partial flags of $\mathbb{F}_{q}-$vector spaces
%$$\{0\}=E^{r}\subset E^{r-1}\subset \ldots \subset E^{1}\subset E^{0}=(\mathbb{F}_{q})^{n}$$ such that ${\rm{dim}}(E^{i-1}/E^{i})=\lambda_{i}$.
%The group $\GL(n,\mathbb{F}_{q})$ acts on $\mathcal{F}_{\lambda}$ in the natural way. %Fix an element $X_{0}\in \mathcal{F}_{\lambda}$ and denote by $\mathcal{P}_{\lambda}$ the stabilizer of $X_{0}$ in $G$ and by $\mathcal{U}_{\lambda}$ the subgroup of elements $g\in \mathcal{P}_{\lambda}$ which induces the identity on $E^{i}/E^{i+1}$ for all $i=0, 1,\ldots, n-1$. Put $\mathcal{L}_{\lambda}=\GL_{\lambda_{r}}(\mathbb{F}_{q})\times \ldots, \times \GL_{\lambda_{1}}(\mathbb{F}_{q})$, then we have $\mathcal{P}_{\lambda}=\mathcal{L}_{\lambda}\times \mathcal{U}_{\lambda}$.

The adjacency matrix of the graph is interpreted as the incidence matrix of the design. Recall that the adjacency matrix $A$ of a multigraph is a $n\times n$ matrix (where $n=|V|$) with rows and columns indexed by the elements of the vertex set and the $(x,y)$- entry is the number of edges connecting $x$ and $y$. If the graph is directed, the matrix $A$ is symmetric and therefore all its eigenvalues are real. The degree of a vertex $deg(v)$ is the number of edges incident with $v$, where we count a loop with multiplicity 2. The largest eigenvalue $\lambda$ of the adjacency matrix describes the spectrum character of the graph topology. 

Given a $t$-design you can associate to it a regular graph, where the points are the nodes of the graph, all the nodes have the same degree and two different nodes are connected if and only if they are in the same block of the design, that is, the neighbors of the vertices are the blocks. Reciprocally given a $k$-regular graph on $v$ vertices, if any two distinct vertices have exactly $\lambda$ common neighbors it is a 2-$(v,k,\lambda)$ design. % If any two adjacent vertices are together adjacent to $\lambda$ vertices, while any two non-adjacent vertices are together adjacent to $\mu$ vertices, is a divisible design. 
It is also possible to design a code which matches the network graph. Then ${\rm{Aut}}(G)$ coincides with the  automorphism group ${\rm{Aut}}(D)$ of the design.

\subsection{Set systems}
A set system is a pair $(X,\mathcal{A})$ such that $X$ is a finite set of points and $\mathcal{A}$ is a set of subsets of $X$, called blocks. The number of points, $|X|$, is the order of the set system. Let $K$ be a set of positive integers. A set system $(X,\mathcal{A})$ is said to be $K$-uniform if $|A|\in K$ for all $A\in \mathcal{A}$. Let $\mathcal{G}=\{G_{1},\ldots, G_{s}\}$ be a partition of $X$ into subsets called groups. The triple $(X,\mathcal{G},\mathcal{A})$ is a group divisible design $(GDD)$ when every 2-subset of $X$ not contained in a group appear in exactly one block and $|A\bigcup\,G|\leq 1$ for all $A\in \mathcal{A}$ and $G\in \mathcal{G}$. A 3-GDD in which all the groups are of size 1 is known as a Steiner triple system.% of order $n$.

\begin{prop} There is bijective correspondence between ordered basis sets of $(\mathbb{F}_{q})^{n}$ and set systems of order $n$. \end{prop}
{\it Proof.} 
This correspondence can be established by associating to any list of $t$ elements contained in $\GL(n,q)$ a partition of $t$ groups of size the order of the corresponding element in $\GL(n,q)$. Namely, to any list $\{\gamma_{1},\ldots, \gamma{t}\}$ of $t$ elements we associate the subgroup $G_{\lambda}$ generated by these $t$ elements. This is a group of type $\lambda$ the partition of orders $\lambda_{i}=ord(\gamma_{i})$ ordered in increasing order $\lambda_{1}\geq \lambda_{2}\geq\ldots \lambda_{t}> 0$.% and $m=\lambda_{1}+\lambda_{2}+ \ldots +\lambda_{t}$.
We assume that $n\geq q-1$ and $G$ is a group containing a Singer cycle $\alpha\in \GL(n,q)$. Let $\Gamma(G_{\lambda})$ be the Cayley graph attached to the subgroup $G_{\lambda}$, that is, the graph in which vertices 1 through $t$ corresponding to each generator are placed in a row with each vertex connected by an unlabelled edge of its immediate neighbors . There is an action of the symmetric group $S_{n}$ on the combinatorial class $\mathcal{G}_{n}$ of regular graphs with $n$ vertices. For any $\sigma\in S_{n}$ and $g\in \mathcal{G}_{n}$, the graph $\sigma\cdot g$ has the same vertex set and edge set as $g$, but each label $i$ in $g$ is replaced by $\sigma^{-1}(i)$ in $\sigma\cdot g$, they are isomorphic graphs.
We define the following linear map over $(\mathbb{F}_{q})^{n}$:
\begin{equation}\label{eq1}\Phi(\Gamma(G_{\lambda}))(x)=A^{G_{\lambda}}_{t,k}\,x. \end{equation}
Here $A^{G}_{t,k}$ is the adjacency matrix of graph $\Gamma(G)$, thus it is a $\{0,1\}$ matrix with rows and columns indexed by the $t$-subspaces and the $k$-subspaces of $\mathbb{F}^{n}_{q}$. 
In particular, constructing $t-$designs over $\mathbb{F}_{q}$ is equivalent to solving the systems of linear Diophantine equations \ref{eq1}. There is a 1 in row $X$ and column $Y$ of M iff $t-$subspaces $X$ is contained in $k-$subspaces Y. With this definition, a $t-(n,k,\lambda)$ design over $\mathbb{F}_{q}$ is precisely a $\{0,1\}$ solution to $A^{G}_{t,k}x=(\lambda, \lambda, \ldots, \lambda)^{T}$, where $\lambda$ is the number of $k-$subspaces containing at least a $t$-subspace, in particular $rank\,(A^{G}_{t,k})\geq t$. 

\cqd

%\begin{example}% From a graph theoretical point of view, 
\subsection{r-designs constructed from the projective line}

Let $X$ be a $v$-set and $\mathcal{P}_{k}(X)$ denote the set of all $k$-subsets of $X$. 
A $t-(v,k,\lambda)$-design is a set system $\mathcal{D}=(X,D)$ in which $D$ is a collection of $\mathcal{P}_{k}(X)$ (called blocks) such that every $t$-subset of $X$ appears in exactly $\lambda$-blocks
A 2-$(v,k,\lambda)$ design is a collection $\mathcal{B}$ of elements of $\mathcal{P}_{k}(X)$ (called blocks) such that every line of the incidence structure $(\mathcal{P}(X),\mathcal{B}(X), I)$ intersect $\mathcal{B}$ in exactly $\lambda$ points.
A 3-$(v,k,\lambda)$ design is a collection of $\mathcal{B}$ of elements of $\mathcal{P}_{k}(X)$ (called blocks) such that any triple $(r_{1},r_{2},r_{3})$ of points is collinear. Such sets are called lines in $\mathcal{B}$ and every line intersect $\mathcal{B}$ in exactly $\lambda$ points. In general $r$-designs admitting $\PG(2,q)$ as a group of automorphisms are known as $(k;r)$ arcs.%we talk about $r-$designs.

% Denote by $l=l(p,q)$ the set of all such points, then $l^{3}\in \mathcal{B}$, that is, any triple $(r_{1},r_{2},r_{3})$ of points in $l$ is collinear. Such sets are called lines in $\mathcal{B}$.

Let $V$ be a 3-dimensional vector space over $\mathbb{F}_{q}$ and consider the projective plane $PG(2,q)$ defined by the incidence structure $(\mathcal{P}(V),\mathcal{B}(V), I)$.
%\end{example}

%Then ${\rm{Aut}}(G)$ coincides with the  automorphism group ${\rm{Aut}}(D)$ of the design.
\begin{defi} A $(k;r)$-arc $\mathcal{K}$ in $PG(2,q)$ is a set of $k$-points such that some $r$, but not $r+1$ of them are collinear. In other words, some line of the plane meets $\mathcal{K}$ in $r$ points and no more than $r$-points. A $(k;r)$-arc is complete if there is no $(k+1;r)$ arc containing it.
\end{defi}
%The $(k;r)-$arcs are a generalization of plane sections for $k-$tuples of points in $\mathcal{P}$ 

\begin{defi} A $k$-arc is a $(k;n,n-1;n,p)$ set with $n\geq 3$ of $k$-points such that, every subset of $s$ points with $s\leq n$ points is linearly independent.
\end{defi}
%If $p\geq n+2$, the NRC is an example of a $(p+1)-$arc. It contains $p+1$ points, and every set of $n+1$ points are linearly independent.

Following the classification of conjugacy classes in $PG(2,q)$ in \cite{SG}, next Lemma classifies designs constructed from the projective line.
\begin{lemma} There are 3 types of $r$-%$(q+1, k,\lambda)$ 
designs constructed from $\PG(2,q)$: unipotent type, semisimple split or semisimple non-split according to the eigenvalues of the representation matrix of the generating elements in $\PG(2,q)$.% into unipotent type if the Jordan form of the representation matrix admits a unique eigenvalue, of semisimple split if the Jordan form of the representation matrix admits 2 eigenvalues o
\end{lemma}
{\it Proof.} If the characteristic polynomial $P(\lambda)$ of the representation matrix $A$ has only one root, call it $\alpha$, it is a primitive element of order $p$ a prime number, then the derived design is called unipotent. It is an arc containing $p+1$ points and for $n<p$ every set of $n+1$ points are linearly independent. If $P(\lambda)$ has two different roots $a, a^{-1}\in \mathbb{F}^{*}_{q}$, $tr(A)=a+a^{-1}$ is an element $\alpha$ of order dividing $\frac{q-1}{d}$. The corresponding design is called semisimple split, and finally if there are no roots, $tr(A)=a+a^{q}=\alpha$, where $a\in \mathbb{F}^{*}_{q^{2}}\backslash \mathbb{F}^{*}_{q}$ is an element $\alpha$ dividing $\frac{q+1}{d}$, the corresponding design is called semisimple non-split.
\cqd

We associate to the 2-design generated by $\tau$ and $\sigma$ the graph which has as vertex set $V$ the points of the projective system $\mathbb{P}((\mathbb{F}_{2})^{m})$ and edge set $E\subseteq [V]^{2}$ the lines of the projective space which corresponds to the blocks of the design. There are $\left[
\begin{matrix}
 m \\
 2\\
\end{matrix}
\right]_{q}$ lines. For any two points there are as much blocks (lines) containing these points as eigenspaces $W_{j}$  by the action of the linear operators $\tau$ and $\sigma$. This special design with parameters $t=2$ and $k=3$ is a Steiner triple system.
The automorphism group of the projective line $\mathbb{P}(\mathbb{F}_{q})$ is the projective linear group $\PGL(2,q)$. Any finite subgroup $A\subset \PGL(2,q)$ defines a $k-$uniform Cayley (sum) hypergraph $\Gamma^{k}(A)$ whose vertices are the generating $k-$tuples of $A$ and the edges are $k-$element sets $\{x_{1},\ldots,x_{k}\}\in {G\choose k} $ represented by random variables $x_{1},\ldots, x_{k}$. In particular, if $f(z)$ is the ordinary generating function that enumerates $A$, that is, number of conjugacy  classes in $A$, then $\frac{1}{1-f(z)}$ is the ordinary generating function enumerating sequences of $k$ elements in $A$. If $G$ is an abelian group, then $x_{1}+\cdots+x_{k}\in A$. In general, we will consider $k$-arcs in $\Gamma(A)$ which represent casual connections between the variables.

The group $\GL(n,q)$ acts transitively on subsets of size $n+1$ of the projective line whenever $q\equiv n+1 (mod\, n+2)$. We can construct secret sharing schemes from configurations of points of size $n+1$ on the projective line. Moreover we can construct secret sharing schemes from configuration of points on curves admitting a transitive linear action. Let $p$ be a prime number and $p\geq n+2$, then the Normal Rational Curve defined as:
%The normal rational curve is defined as:
$$\mathcal{V}^{n}_{1}:=\Big\{F (1,x,x^{2},\ldots,x^{n})| \ x\in \mathbb{F}_{p}\bigcup \{\infty\}\Big\}$$ 
is an example of a $(p+1)$-arc. It contains $p+1$ points, and every set of $n+1$ points are linearly independent.

\end{document}